\documentclass[12equipped]{article}
\usepackage{amssymb,amsmath}

\input epsf.tex

\title{Quantum p-adic spaces and quantum p-adic groups}

\author {Yan Soibelman}

\vspace{3mm}

\begin{document}
\maketitle
\newtheorem{thm}{Theorem}[subsection]
\newtheorem{defn}[thm]{Definition}
\newtheorem{lmm}[thm]{Lemma}
\newtheorem{rmk}[thm]{Remark}
\newtheorem{prp}[thm]{Proposition}
\newtheorem{conj}[thm]{Conjecture}
\newtheorem{exa}[thm]{Example}
\newtheorem{cor}[thm]{Corollary}
\newtheorem{que}[thm]{Question}
\newtheorem{ack}{Acknowledgements}
\newcommand{\C}{{\bf C}}
\newcommand{\K}{{\bf k}}
\newcommand{\R}{{\bf R}}
\newcommand{\N}{{\bf N}}
\newcommand{\Z}{{\bf Z}}
\newcommand{\Q}{{\bf Q}}
\newcommand{\G}{\Gamma}
\newcommand{\A}{A_{\infty}}
\newcommand{\ihom}{\underline{Hom}}
\newcommand{\epi}{\twoheadrightarrow}
\newcommand{\mono}{\hookrightarrow}
\newcommand\ra{\rightarrow}
\renewcommand\O{{\cal O}}
\newcommand\nca{nc{\bf A}^{0|1}}
\newcommand{\epp}{\varepsilon}

{\it To the memory of Sasha Reznikov}

\tableofcontents

\section{Introduction}

The paper is devoted to examples of quantum spaces over non-archimedean fields and is, in a sense,
a continuation of [So1] (part of the material is borrowed from the loc. cit).
There are three classes of examples which I discuss in this paper: quantum affinoid spaces, quantum non-archimedean Calabi-Yau varieties
and quantum $p$-adic groups. Let us recall the definitions and discuss the contents of the paper.

Quantum affinoid algebras  are defined similarly to
the  ``classical case" $q=1$. It is a special case of a more general notion of non-commutative affinoid
algebra introduced in [So1]. Let $k$ be a Banach field and
${k}\langle\langle T_1,...,T_n\rangle\rangle$ be the algebra
of formal series in free variables $T_1,...,T_n$.
For each $r=(r_1,...,r_n),r_i\ge 0, 1\le i\le n$
we define a subspace ${k}\langle\langle T_1,...,T_n\rangle\rangle_r$
consisting of series $f=\sum_{i_1,...,i_m}a_{i_1,...,i_m}T_{i_1}...T_{i_m}$
such that $\sum_{i_1,...,i_m}|a_{i_1,...,i_m}|r_{i_1}...r_{i_m}<+\infty$.
Here the summation is taken over all sequences $(i_1,...,i_m), m\ge 0$ and
$|\bullet|$ denotes the norm in $k$.
In this paper we consider the case when 
$k$ is a valuation field (i.e. a Banach field with respect to a  multiplicative non-archimedean norm). 
In the non-archimedean case the convergency condition is replaced
by the following one: $max\,|a_{i_1,...,i_m}|r_{i_1}...r_{i_m}\to 0$ as $i_1+...+i_m\to \infty$.
Clearly each ${k}\langle\langle x_1,...,x_n\rangle\rangle_r$ is a Banach
algebra called the algebra of analytic functions on a non-commutative
$k$-polydisc $E_{NC}(0,r)$ centered at zero and having the (multi)radius $r=(r_1,...,r_n)$.
The norm is given by $max\,|a_{i_1,...,i_m}|r_{i_1}...r_{i_m}$.
A {\it non-commutative $k$-affinoid algebra} is an admissible noetherian quotient of this algebra (cf. [Be1],
Definition 2.1.1).
Let us fix $q\in k^{\ast}$ such that $|q|=1$, and  $r=(r_1,...,r_n), r_i\ge 0$.
A {\it quantum $k$-affinoid algebra} is a special case of the previous definition. It is
defined as an admissible quotient 
of the algebra $k\{T\}_{q,r}:=k\{T_1,...,T_n\}_{q,r}$ of the series $f=\sum_{l=(l_1,...,l_n)\in \Z_+^n}a_lT_1^{l_1}...T_n^{l_n}$
such that $a_l\in k$, $T_iT_j=qT_jT_i, j<i$, and $max|a_l|r^{|l|}\to 0$ as $|l|:=l_1+...+l_n\to \infty$.
The latter is also called the algebra of analytic functions on the quantum polydisc 
$E_{q}(0,r)$. It is less useful notion than the one of non-commutative affinoid algebra since there are few two-sided closed ideals in the algebra
$k\{T_1,...,T_n\}_{q,r}$. Nevertheless quantum affinoid algebras appear in practice (e.g. in the case
of quantum Calabi-Yau manifolds considered below).
In the case when all $r_i=1$ we speak about {\it strictly $k$-affinoid}
non-commutative (resp. quantum) algebras, similarly to [Be1]. Any non-archimedean extension $K$ of $k$ gives rise to a
non-commutative (resp. quantum) affinoid $k$-algebra, cf. loc.cit. 
There is a generalization of quantum affinoid algebras which we will also call quantum affinoid algebras.
Namely, let $Q=((q_{ij}))$ be an $n\times n$ matrix with
entries from $k$ such that $q_{ij}q_{ji}=1, |q_{ij}|=1$ for all $i,j$.
Then we define the quantum affinoid algebra as an admissible quotient of the algebra
$k\{T_1,...,T_n\}_{Q,r}$. The latter defined similarly to $k\{T_1,...,T_n\}_{q,r}$, but now we use
polynomials in variables $T_i, 1\le i\le n$
such that $T_iT_j=q_{ij}T_jT_i$.
One can  think of $k\{T_1,...,T_n\}_{Q,r}$ as of the
quotient of
$k\langle\langle T_i,t_{ij}\rangle \rangle_{r,{\bf 1}_{ij}}$,
where $1\le i,j\le n$ and ${\bf 1}_{ij}$ is the unit $n\times n$
matrix,
by the  two-sided ideal generated by the relations
$$t_{ij}t_{ji}=1, T_iT_j=t_{ij}T_jT_i, t_{ij}a=at_{ij},$$
for all indices $i,j$ and all
$a\in k\langle\langle T_i,t_{ij}\rangle \rangle_{r,{\bf 1}_{ij}}$.
In other words, we treat $q_{ij}$ as variables
which belong to the center of our algebra and have the
norms equal to one. Having the above-discussed generalizations of affinoid algebras
we can consider their Berkovich spectra (sets of multiplicative seminorms).
Differently from the commutative case, the theory
of non-commutative and quantum analytic spaces is not developed yet (see discussion in [So1]).
 
Quantum Calabi-Yau manifolds provide examples of topological spaces equipped with rings of
non-commutative affinoid algebras, which are quantum affinoid outside of a ``small" subspace.
More precisely,
let  $k=\C((t))$ be the field of Laurent series, equipped with its standard valuation (order
of the pole) and the corresponding non-archimedean norm.
Quantum Calabi-Yau manifold of dimension $n$ over $\C((t))$, is defined as a ringed space $(X,{\cal O}_{q,X})$ which
consists of an analytic Calabi-Yau manifold $X$ of dimension $n$ over $\C((t))$ and a sheaf of
$\C((t))$-algebras ${\cal O}_{q,X}$ on $X$
such that ${\cal O}_{q,X}(U)$ is a non-commutative affinoid algebra for any affinoid $U\subset X$
and the following two conditions are satisfied:

1) $X$ is a $\C((t))$-analytic manifold corresponding to a maximally degenerate algebraic Calabi-Yau manifold
$X^{alg}$  of dimension $n$ (see [KoSo2] for the definitions);

2) Let $Sk(X)$ be the skeleton of $X$ defined in [KoSo1], and let us choose a projection $\pi: X\to Sk(X)$ 
described in the loc.cit. Then the direct image $\pi_{\ast}({\cal O}_{q,X})$ is locally isomorphic
(outside of a topological subvariety of the codimension at least two) to the 
sheaf of $\C((t))$-algebras ${\cal O}_{q,{\R}^n}^{can}$ on ${\R}^n$ which is characterized by the property
that for any open connected subset $U\subset {\R}^n$ we have
${\cal O}_{q,{\R}^n}^{can}(U)=\{\sum_{l\in {\Z}^n}a_lz^l\}$ such that $a_l\in \C((t))$ and
$sup_{l\in{\Z}^n} (log|a_l|+\langle l,x\rangle)<\infty$ for any $x\in U$. Here $\langle(l_1,...,l_n),(x_1,...,x_n)\rangle=\sum_{1\le i\le n}l_ix_i$;

For a motivation of this definition see [KoSo1-2] in the ``commutative" case $q=1$. Roughly speaking, in that case 
the above definition requires the Calabi-Yau manifold $X$ to be locally isomorphic (outside
of a ``small" subspace) to an analytic torus fibration $\pi_{can}:({\bf G}_m^{an})^n\to {\R}^n$,
where on $\C((t))$-points the canonical projection $\pi_{can}$ is the ``tropical" map
$(z_1,...,z_n)\mapsto (log|z_1|,..., log|z_n|)$. This is a ``rigid-analytic" implementation of the Strominger-Yau-Zaslow
conjecture in Mirror Symmetry (see [KoSo1-2] for more on this topic).
In present paper we discuss the case $n=2$, essentially following [KoSo1], [So1]. 
Perhaps the higher-dimensional  case
can be studied by the technique developed in a recent paper [GroSie1] (which  in some sense
generalizes to the higher-dimensional case ideas of [KoSo1]). We plan to return to this problem
in the future.

Finally, we discuss the notion of $p$-adic quantum group.
Quantum groups over  $p$-adic fields and their representations will be discussed in more detail in the forthcoming paper [So2].
We have borrowed some material from there. Recall that quantum groups are considered
in the literature either in the framework of algebraic groups or in some special examples of locally compact groups over $\R$. In the  case of groups over $\R$ or $\C$ there is the following problem:
how to describe, say, smooth or analytic (or rapidly descreasing) functions on a complex or real Lie group in terms of the representation
theory of its enveloping algebra? Finite-dimensional representations give rise to the algebra of regular functions (via Peter-Weyl theorem), but more general classes of functions are not so easy to handle.
The case of $p$-adic fields is different for two reasons. First, choosing a good basis in the
enveloping algebra, we can consider series with certain restrictions on the growth of norms
of their coefficients. This allows us to describe a basis of compact open neighbourhoods of the unit of
the corresponding $p$-adic group. Furthermore, combining the ideas of [ShT1] with the approach of [So3]
one can define the algebra of {\it locally analytic}  functions on a compact $p$-adic group as a certain completion
of the coordinate ring of the group. Dualizing, one obtains the algebra of locally-analytic distributions.
According to [ShT1] modules over the latter provides an interesting class of $p$-adic representations,
which contains, e.g. principal series representations. The above considerations can be ``quantized",
giving rise to quantum locally-analytic groups.

Present paper contains a discussion of the above-mentioned three classes of examples of non-commutative spaces.
The proofs are omitted and will appear in separate publications. We should warn the reader that the paper
does not present a piece of developed theory. This explains its sketchy character.
My aim is to show interesting classes of non-archimedean non-commutative spaces which can be obtained
as analytic non-commutative deformations of the corresponding classical spaces. 
They deserve  further study (for the quantum groups case see [So2]).

When talking about rigid analytic spaces we use the approach of Berkovich, which seems to be more suitable in
the non-commutative framework. For this reason our terminology is consistent with [Be1].

{\it Acknowledgements.} I am grateful to many people who shared with me their ideas and insights,
especially to Vladimir Berkovich, Joseph Bernstein, Matthew Emerton and Maxim Kontsevich.
Excellent papers [SchT1-6] by Peter Schneider and Jeremy Teitelbaum played a crusial role in convincing
me that the theory of quantum p-adic groups should exist. 
I thank IHES for the hospitality and excellent research conditions. This work was partially supported by an NSF grant.

\section{Quantum affinoid algebras}

Let $k$ be a  valuation field.

We recall here the definition already given in the Introduction.
Let us fix $r=(r_1,...,r_n)\in \R_{\ge 0}^{n}$.
We start with the algebra
$k\langle T\rangle:=k\langle T_1,...,T_n\rangle$ of polynomials
in $n$ free variables
and consider its completion $k\langle \langle T\rangle\rangle_r$ with respect to the norm
$|\sum_{\lambda\in P({\Z}_+^n)}a_{\lambda}T^{\lambda}|=max_{\lambda} |a_{\lambda}|r^{\lambda}$.
Here $P({\Z}_+^n)$ is the set of finite paths in ${\Z}_+^n$ starting
at the origin, and $T^{\lambda}=T_1^{\lambda_1}T_2^{\lambda_2}....$ for the path
which moves $\lambda_1$ steps in the direction $(1,0,0,...)$ then $\lambda_2$
steps in the direction $(0,1,0,0...)$, and so on (repetitions are allowed, so we
can have a monomial like $T_1^{\lambda_1}T_2^{\lambda_2}T_1^{\lambda_3}$).

\begin{defn}
We say that a noetherian Banach unital algebra $A$ is {\it non-commutative affinoid $k$-algebra}
if there is an admissible surjective homomorphism
$k\langle \langle T\rangle\rangle_r\to A$ (admissibility
means that the norm on the image is the quotient norm).

\end{defn}

In particular, affinoid algebras in the sense of [Be1] belong to this class
(unfortunately the terminology is confusing in this case: commutative affinoid
algebras give examples of non-commutative affinoid algebras!).
Another class of examples is formed by quantum affinoid algebras defined in the Introduction.

Let us now recall the following definition (see [Be1]).

\begin{defn} Berkovich spectrum $M(A)$ of a unital Banach ring $A$ consists of bounded multiplicative
seminorms on $A$.

\end{defn}

If $A$ is a $k$-algebra, we require that seminorms extend the norm on $k$.
It is well-known (see [Be1], Th. 1.2.1) that if $A$ is commutative then $M(A)$ is a non-empty
compact Hausdorff topological space (in the
weak topology). If $\nu\in M(A)$ then $Ker\,\nu$ is a two-sided closed prime ideal in $A$.
Therefore it is not clear whether $M(A)$ is non-empty in the non-commutative case.

Algebras of analytic functions on the non-commutative and quantum polydiscs
carry multiplicative ``Gauss norms" (see Introduction), hence the Berkovich spectrum
is non-empty in each of those cases.
The following  example  can be found in [So3], [SoVo].

Let $L$ be a free abelian
group of finite rank $d$, $\varphi:L\times L\to {\Z}$ be a skew-symmetric
bilinear form, $q\in K^{\ast}$ satisfies the
condition $|q|=1$. Then $|q^{\varphi(\lambda,\mu)}|=1$
for any $\lambda,\mu \in L$. We denote by
$A_q(T(L,\varphi))$ the
{\it algebra of regular functions on the quantum torus $T_q(L,\varphi)$}.
By definition, it is a $k$-algebra with generators $e(\lambda),\lambda\in L$,
subject to the relation
$$e(\lambda)e(\mu)=q^{ \varphi(\lambda,\mu)}e(\lambda+\mu).$$
The algebra of analytic functions on the analytic quantum torus 
$T_q^{an}(L,\varphi)$ consists by definition of series 
$\sum_{\lambda\in L}a(\lambda)e(\lambda)$,$a(\lambda)\in k$
such that for all $r=(r_1,...,r_d), r_i>0$ one has:
$|a(\lambda)|r^{\lambda}\to 0$ as $|\lambda|\to \infty$ (here $|(\lambda_1,...,\lambda_d)|=
\sum_i|\lambda_i|$).

Quantum affinoid algebra $k\{T\}_{q,r}$ discussed in the Introduction is the algebra
of analytic functions on quantum polydisc of the (multi)radius $r=(r_1,...,r_n)$.
It was shown in [So1] that $M(k\{T\}_{q,r})$ can be quite big as long as $|q-1|<1$.
In particular, it contains ``quantum" analogs of the norms $|f|_{E(a,\rho)}$
which is the ``maximum norm" of an analytic function $f$ on the polydisc centered at
$a=(a_1,...,a_n)$ of the radius $\rho=(\rho_1,...,\rho_n)$, with the condition
$a_i\le \rho_i<r_i, 1\le i\le n$. Similar result holds for the quantum analytic torus.
This observation demonstrates an interesting phenomenon: differently from the formal deformation
quantization, the non-archimedean analytic quantization ``preserves" some part
of the spectrum of the ``classical" object.

The conventional definition of the quatization can be carried out to the analytic case with obvious
changes. Indeed, the notion of Poisson algebra admits a straightforward generalization to the analytic case
(Poisson bracket is required to be a  bi-analytic map). Furthermore, for any commutative affinoid algebra
$A$ there is a notion of non-commutative $A$-affinoid algebra, which is a natural generalization
of the notion of $k$-affinoid algebra (we use $A\langle\langle T_1,...,T_n\rangle\rangle_r$ instead
of $k\langle\langle T_1,...,T_n\rangle\rangle_r$).

Let now ${\cal O}(E(0,r))$ be the algebra of analytic functions on a $1$-dimensional polydisc 
$E(0,r)=M(k\{r^{-1}T\})$ 
of the radius $r$ (the notation is from [Be1], Chapter 2). We say that a non-commutative 
${\cal O}(E(0,r))$-affinoid algebra $A$ is an {\it 
analytic quantization of a $k$-affinoid commutative Poisson
algebra $A_0$ over the polydisc $E(0,r)$} if the following two conditions are satisfied:

1) $A$ is a topological ${\cal O}(E(0,r))$-algebra, free as a topological ${\cal O}(E(0,r))$-module.

2) The quotient $A/TA$ is isomorphic to $A_0$ as a $k$-affinoid Poisson algebra.

Then a quantization of a $k$-analytic space $(X,{\cal O}_X)$
iss a ringed space $(X,{\cal O}_{q,X})$ such that for any affinoid $U\subset X$ the
algebra ${\cal O}_{q,X}(U)$ is an analytic quantization of ${\cal O}_{X}(U)$ over some polydisc
$E(0,r)$.

Notice that the projection $A\to A_0$ induces an embedding of Berkovich spectra
$M(A_0)\to M(A)$. Every element of $A$ can be thought of as analytic function
on $E(0,r)$ with values in a non-commutative $k$-affinoid algebra.
Suppose that $A\simeq A_0\{r^{-1}T\}$ as a $k\{r^{-1}T\}$-module.
Then the topological vector space $A$ is isomorphic to the space of analytic functions
on $E(0,r)$ with values in $A_0$ (but the product is not a pointwise product of functions).
Assume that $r\le 1$ and consider the subspace $A_1$ of analytic functions $a(x)$ as above such that
$|a(0)|_{A_0}\le 1, |a(x)-a(0)|_{A_0}\le |T(x)|, x\in E(0,r)$, where $\bullet|_{A_0}$ denotes
the norm on $A_0$. Here $x$ is interpreted as a seminorm on the
Banach $k$-algebra   $k\{r^{-1}T\}$, hence $|T(x)|$ is the norm of the generator $T$ in the
completition of the residue
field $k\{r^{-1}T\}/Ker\,x$. It is clear that $A_1$ is in fact a Banach $k$-algebra.
Hence the natural projection $a(x)\mapsto a(0)$ defines an embedding
$M(A_0)\to M(A_1)$.

Suppose that $X$ is an analytic spaces for which there is a notion of a skeleton $Sk(X)$ either
in the sense of [KoSo1] (then $X$ is assumed to be Calabi-Yau) or in the sense of [Be2,Be3].
Then in either of these cases there is a continuous retraction $\pi: X\to Sk(X)$.
Suppose that the there is a quantization $(X,{\cal O}_{q,X})$ of $(X,{\cal O}_X)$ in the above sense.

\begin{conj} For any closed $V\subset X$ there is a natural embedding 
$i_V:V\subset M({\cal O}_{q,X}(\pi^{-1}(V)))$ such that $\pi\circ i_V=id_V$.
Moreover if $V_1\subset V_2$ then the restriction of $i_{V_1}$ to $V_2$ is equal
to $V_2$.

\end{conj}

In other words, the skeleton survives an analytic quantization.
The above conjecture is not very precise, because there is no general definition
of a skeleton. The definition given in [KoSo1] is different from the one in [Be2,3] even for Calabi-Yau manifolds. Hence the conjecture is an ``experimental fact" at this time. 

\section{Quantum Calabi-Yau varieties}

Let $X^{alg}$ be a maximally degenerate (in the sense of [KoSo1-2]) algebraic Calabi-Yau  manifold
over $\C((t))$ of dimension $n$ and $X$ be the corresponding $\C((t))$-analytic space. 
Then one can associate with $X$ a $PL$-manifold
$Sk(X)$ of real dimension $n$, called the skeleton of $X$ (see [KoSo1-2]).
A choice of K\"ahler structure on $X^{alg}$ defines (conjecturally) a continuous retraction
$\pi: X\to Sk(X)$.
This map satisfies the condition 2) from the Introduction. In other words, it defines
a (singular) analytic torus fibration over $Sk(X)$ with the generic fiber, isomorphic
to the analytic space $M(k[T_1^{\pm 1},...,T_n^{\pm 1}]/(|T_i|=c_i, 1\le i\le n))$, 
where $c_i>0, 1\le i\le n$ are some numbers.
Since the projection is Stein (see loc. cit), one can reconstruct $X$ (as a ringed space)
from the knowledge of $(Sk(X),\pi_{\ast}({\cal O}_{X}))$, where 
${\cal O}_{X}$ is the sheaf of analytic functions on $X$.

Let $B=Sk(X)$ and $B^{sing}$ be the ``singular subvariety" of real codimension two (see Introduction).
It was observed in [KoSo1] that the norms of elements
of the direct image sheaf $\pi_{\ast}({\cal O}_{X}^{\times})$ define an integral affine structure on 
$B^0:=B\setminus B^{sing}$.
Hence we would like to reconstruct the analytic space starting with a PL-manifold equipped with
a (singular) integral affine structure.  As we will explain in the next subsection the same data
give rise to a sheaf of quantum affinoid algebras on $B^0$.

\subsection{Integral affine structures and quantized canonical sheaf}

Here we explain following [KoSo1] and [So1] how a manifold with integral affine structure
defines a sheaf of (quantum) affinoid algebras.

Recall that an integral affine structure ($\Z$-affine structure
for short)
on an $n$-dimensional topological manifold $Y$ is given by
a maximal atlas of  charts such that the change of coordinates
between any two charts is described by the formula
$$ x_i^{\prime}=\sum_{1\le j\le n}a_{ij}x_j+b_i,$$
where $(a_{ij})\in GL(n,{\Z}), (b_i)\in {\R}^n$.
In this case one can speak about the sheaf of $\Z$-affine functions,
i.e. those which can be locally  expressed in affine coordinates by the formula
$f=\sum_{1\le i\le n}a_ix_i+b, a_i\in \Z, b\in \R.$
An equivalent description: $\Z$-affine structure is given by
a covariant lattice $T^{\Z}\subset TY$ in the tangent bundle (recall that an affine
structure on $Y$ is the same as a torsion free flat connection on the
tangent bundle $TY$).

Let $Y$ be a manifold with $\Z$-affine
structure. The sheaf of $\Z$-affine functions
$Aff_{\Z}:=Aff_{{\Z},Y}$ gives rise to an exact
sequence of sheaves of abelian groups

$$ 0\to {\R}\to Aff_{\Z}\to (T^{\ast})^{\Z}\to 0,$$

where $(T^{\ast})^{\Z}$ is the sheaf associated with the dual to the covariant lattice
$T^{\Z}\subset TY$.

Let us recall the following notion introduced in [KoSo1], Section 7.1. Let $k$ be a valuation field.

\begin{defn} A  $k$-affine structure on $Y$
compatible with the given $\Z$-affine structure
is a sheaf $Aff_k$ of abelian groups
on $Y$, an exact sequence of sheaves

$$1\to k^{\times}\to Aff_k\to (T^{\ast})^{\Z}\to 1,$$
together with a homomorphism $\Phi$ of this exact sequence
to the exact sequence of sheaves of abelian groups
$$0\to {\R}\to Aff_{\Z}\to (T^{\ast})^{\Z}\to 0,$$
such that  $\Phi=id$ on $(T^{*})^{\Z}$ and
$\Phi=val$ on $k^{\times}$, where $val$ is the valuation map.

\end{defn}

Since $Y$ carries a $\Z$-affine structure, we
have the corresponding $GL(n,{\Z})\ltimes {\R}^n$-torsor on $Y$,
whose fiber
over a point $x$ consists of all $\Z$-affine coordinate systems
at $x$.

Then one has the following equivalent description of the notion
of $k$-affine structure.

\begin{defn} A $k$-affine structure on $Y$
compatible with the given $\Z$-affine structure is
a $GL(n,{\Z})\ltimes (k^{\times})^n$-torsor
on $Y$ such that the application of $val^{\times n}$ to
$(k^{\times})^n$ gives the initial $GL(n,{\Z})\ltimes {\R}^n$-torsor.
\end{defn}

Assume that $Y$ is oriented and carries a $k$-affine structure compatible
with a given $\Z$-affine structure.
Orientation allows us to reduce
the structure group of the torsor
defining the $k$-affine structure to $SL(n,{\Z})\ltimes (k^{\times})^n$ .

Let $q\in k, |q|=1$, and $z_1,...,z_n$ be invertible variables such that $z_iz_j=qz_jz_i$,
for all $1\le i<j\le n$. We define the sheaf of $k$-algebras ${\cal O}^{can}_q$
on ${\R}^n, n\ge 2$ by the formulas:
$${\cal O}^{can}_q(U)=\left\{
\sum_{I=(I_1,...,I_n)\in {\Z}^n}c_{I}z^I,|\,\forall (x_1,...,x_n)\in U\,\,\,\sup_{I}\left(
\log(|c_{I}|)+\sum_{1\le m\le n}I_mx_m\right)<\infty\right\},$$
where $z^{I}=z_1^{I_1}\dots z_n^{I_n}.$
Since $|q|=1$ the convergency
condition does not depend on the order of variables.

The sheaf ${\cal O}^{can}_q$ can be lifted to $Y$ (we keep the same
notation for the lifting).
In order to do that it suffices to define the action of the group $SL(n,{\Z})\ltimes (k^{\times})^n$
on the canonical sheaf on ${\R}^n$.
Namely, the inverse to an element
$(A,\lambda_1,...,\lambda_n)\in SL(n,{\Z})\ltimes (k^{\times})^n$
acts on monomials as
$$z^{I}=z_1^{I_1}\dots z_n^{I_n}\mapsto \left({\textstyle\prod_{i=1}^n}\lambda_i^{I_i}\right)\,\,z^{A(I)}\,\,\,.$$
The action
of the same element on ${\R}^n$ is given by a similar formula:
$$x=(x_1,\dots ,x_n)\mapsto A(x)-(val(\lambda_1),\dots,val(\lambda_n))\,\,\,.$$

Any $n$-dimensional manifold $Y$ with integral affine structure admits a covering by charts
with transition functions being integral affine transformations. This allows to define the sheaf
${\cal O}_{q,Y}^{can}$ as the one which is locally isomorphic to 
${\cal O}_q^{can}={\cal O}_{q,{\R}^n}^{can}$.

It is explained in [KoSo1] (see also [So1], Section 7.2) that for any open $U\subset {\R}^n$ 
the topological space $M({\cal O}^{can}_{q}(U))$ for $q=1$ is an analytic torus fibration
in the sense of Introduction. Recall that an analytic torus fibration
is a fiber bundle $(X,Y,\pi)$ consisting of a commutative $k$-analytic space,
a topological manifold $Y$ and a continuous map $\pi: X\to Y$ such that it is locally
isomorphic to the torus fibration ${\bf G}_m^n\to {\R}^n$ from Introduction.
In that case $\pi$ is a Stein map, and we have: $\pi^{-1}(U)=M({\cal O}^{can}_{q=1}(U))$.
Therefore we can think of the ringed space $(Y,{\cal O}_{q,Y}^{can})$ as of quantization of 
this torus fibration.

\subsection{Model sheaf near a singular point}

In the case of maximally degenerate K3 surfaces the skeleton is homeomorphic
to $B=S^2$ (the two-dimensional sphere) equipped with an integral affine structure outside
of the subset $B^{sing}$ consisting of $24$ points (see [KoSo1], Section 6.4, where the affine structure
is described). The construction of the previous subsection gives rise to a sheaf
of quantum  $\C((t))$-affinoid algebras over $B^0=B\setminus B^{sing}$.
In order to complete the quantization procedure we need to extend the sheaf ${\cal O}_{q,B^0}^{can}$ to a neighbourhood
of $B^{sing}$. It is explained in [KoSo1] (case $q=1$) and in [So1] (case $|q|=1$) that
one has to modify this sheaf in order to extend it to singular points. 
Summarizing,
the quantization is achieved in two steps.
First, we define a sheaf of non-commutative $\C((t))$-affinoid algebras in a neighbourhood of $B^{sing}$ such that
it is locally isomorphic to the canonical sheaf ${\cal O}_{q,B^0}^{can}$ outside of $B^{sing}$ and gives a ``local model"
for the future sheaf $\pi_{\ast}({\cal O}_{q,X})$ at the singularities. 
Second, we modify the sheaf ${\cal O}_{q,B^0}^{can}$ by applying (infinitely many times) automorphisms associated with edges of an
infinite tree embedded in $B$, such that its  external vertices belong to $B^{sing}$. Those modifications ensure
that the resulting sheaf can be glued with the model sheaf at the singularities, and that it is indeed the direct image of the sheaf of analytic functions on a compact $\C((t))$-analytic K3 surface.
More precisely, we do the following.

We start with an open covering of $\R^2$ by the following sets $U_i, 1\le i\le 3$.
Let us fix a number $0<\varepsilon<1$ and define
$$\begin{array} {lll}
U_1 & = & \{(x,y)\in {\R}^2|x<\varepsilon |y|\,\}\\
U_2 & = & \{(x,y)\in {\R}^2|x>0, y<\varepsilon x\,\}\\
U_3 & = & \{(x,y)\in {\R}^2|x>0, y>0\}
\end{array}$$
Clearly ${\R}^2\setminus\{(0,0)\}=U_1\cup U_2\cup U_3$.
We will also need a slightly modified domain $U_2'\subset U_2$ defined as
$\{(x,y)\in {\R}^2|x>0, y<\frac{\varepsilon}{1+\varepsilon} x\,\}$.

Let $\pi_{can}: ({\bf G}_m^{an})^2\to \R^2$ be the canonical map defined in the Introduction
(see also [KoSo1]).
We define the following three open subsets of the two-dimensional analytic torus: $T_i:=\pi_{can}^{-1}(U_i), i=1,3$ and
 $T_2:=\pi_{can}^{-1}(U_2')$.
There are natural projections $\pi_i: T_i\to U_i$  given by the formulas
$$
\pi_i(|\bullet|)=\pi_{can}(|\bullet|)=
(\log|\xi_i|, \log|\eta_i|), \,\,\,i=1,3$$
$$\pi_2(|\bullet|)=\left\{\begin{array}{ll}(\log|\xi_2|, \log|\eta_2|)& \mbox{  if }|\eta_2|<1\\
(\log|\xi_2|-\log|\eta_2|, \log|\eta_2|) & \mbox{ if }
|\eta_2|\ge 1
\end{array}\right. $$

To each $T_i$ we assign the algebra ${\cal O}_q(T_i)$ of series
$\sum_{m,n}c_{mn}\xi_i^m\eta_i^n$ such that $\xi_i\eta_i=q\eta_i\xi_i$, $c_{mn}\in \C((t))$,
and for the seminorm $|\bullet|$ corresponding to a point of $T_i$ (which means
that $(log|\xi_i|,log|\eta_i|)\in U_i$) one has:
$sup_{m,n}(m\,log|\xi_i|+n\,log|\eta_i|)<+\infty$. Similarly, we can define ${\cal O}_q(U)$
for any $U\subset U_i$.
In this way we obtain a sheaf
of quantum $\C((t))$-affinoid  algebras on the set $U_i$. We will denote this sheaf by 
$\pi_{i\ast}({\cal O}_{q,T_i})$.

We define the sheaf ${\cal O}_q^{can}$ on ${\R}^2\setminus \{(0,0)\}$ 
as $\pi_{i*}\left(\O_{q,T_i}\right)$ on each domain $U_i$,
with identifications

$$\begin{array}{llcl}
(\xi_1,\eta_1) & = & (\xi_2,\eta_2) & \mbox{ on } U_1\cap U_2 \\
(\xi_1,\eta_1) & = & (\xi_3,\eta_3) & \mbox{ on } U_1\cap U_3 \\
(\xi_2,\eta_2) & = & (\xi_3\eta_3,\eta_3) & \mbox{ on } U_2\cap U_3
\end{array}$$

The notation for the sheaf is consistent with the previous subsection since
${\cal O}_q^{can}$ is locally isomorphic to the canonical sheaf associated with the standard
integral affine structure.

Let us modify the canonical sheaf ${\cal O}_q^{can}$ in the following way.
On the sets $U_1$ and $U_2\cup U_3$ the new
sheaf ${\cal O}_q^{mod}$ is isomorphic to  ${\cal O}^{can}_q$
(by identifying of coordinates $(\xi_1,\eta_1)$
and glued coordinates $(\xi_2,\eta_2)$ and $(\xi_3,\eta_3)$ respectively).
On the intersection $U_1\cap (U_2\cup U_3)$ we identify two
copies of  the canonical sheaf by an automorphism $\varphi$ of
${\cal O}^{can}_q$ given (we skip the index of the
coordinates) by
 $$\varphi(\xi,\eta)=\left\{\begin{array}{cll}
(\xi(1+\eta),\eta) & \mbox{ on } & U_1\cap U_2 \\
(\xi(1+{\eta}^{-1}),\eta) & \mbox{ on } & U_1\cap U_3
\end{array}
\right.$$

Finally we are going to introduce a sheaf of $\C((t))$-algebras ${\cal O}_q^{sing}$ on 
a small open disc  $W\subset {\R}^2, \{(0,0)\}\in W$ such that
${{\cal O}_q^{sing}}|_{W\setminus \{(0,0)\}}$ is isomorphic to ${{\cal O}_q^{mod}}|_{W\setminus \{(0,0)\}}$.
The sheaf ${\cal O}_q^{sing}$  provides a non-commutative deformation
of the ``local model sheaf" near a singular point (see [KoSo1], Section 8 about the latter).

Let  us consider a non-commutative $\C((t))$-algebra $A_q(S)$ generated
by $\alpha,\beta,\gamma$ subject to the following
relations:
$$\alpha\gamma=q\gamma\alpha,q\beta\gamma=\gamma\beta,$$
$$\beta\alpha-q\alpha\beta=1-q,$$
$$(\alpha\beta-1)\gamma=1.$$
For $q=1$ this algebra coincides with the algebra of regular
functions on the surface $S\subset {\bf A}^3_{\C((t))}$ given
by the equation $(\alpha\beta-1)\gamma=1$ and moreover, it is a flat
deformation of the latter with respect to the parameter $q-1$.
It is explained in [KoSo1], Section 8, that there is a natural map $p: S^{an}\to \R^2$ 
of the corresponding analytic surface such that
$p_{\ast}({\cal O}_{S^{an}})$ is a local model near a singularity of the sheaf $\pi_{\ast}({\cal O}_X)$,
where $X$ is the maximally degenerate K3 surface and $\pi$ is the projection to the skeleton $Sk(X)$. 

Let us denote by ${\cal O}_{q,r_1,r_2,r_3}(S^{an})$ the non-commutative affinoid algebra 
which the quotient of $\C((t))\langle\langle \alpha,\beta,\gamma\rangle\rangle_{r_1,r_2,r_3}$ by the closed two-sided
ideal generated by the above three relations for $A_q(S)$. Here $r_i, i=1,2,3$ are arbitrary non-negative numbers.
We denote by ${\cal O}_{q}(S^{an})$ the intersection of all algebras ${\cal O}_{q,r_1,r_2,r_3}(S^{an})$.

We define homomorphisms of  non-commutative  algebras
$g_i: A_q(S)\to {\cal O}_q(T_i), 1\le i\le 3$ by the following 
formulas (the notation is obvious):

$$\begin{array}{lll}
g_1(\alpha,\beta,\gamma) & = & ({{\xi_1}^{-1}}, \xi_1(1+\eta_1), {\eta_1}^{-1})\\
g_2(\alpha,\beta,\gamma) & = & ({(1+\eta_2){\xi_2}^{-1}}, \xi_2, {\eta_2}^{-1})\\
g_3(\alpha,\beta,\gamma) & = &
((1+\eta_3)(\xi_3\eta_3)^{-1}, \xi_3\eta_3, (\eta_3)^{-1})\end{array}$$

These homomorphisms correspond to the natural embeddings $T_i\mono S^{an}, i=1,2,3$.
One can use these homomorphisms in order to show an existence of non-trivial multiplicative seminorms
on $A_q(S)$ and construct explicitly some of the corresponding representations of $A_q(S)$
in a $k$-Banach vector space.

For example, let us consider a Banach vector space $V_r$ consisting of series
$\sum_{i\in \Z}a_iT^i, a_i\in k$ such that $|a_i|r^i\to 0$ as $|i|\to \infty$, where $r>0$ is some number.
Let $\tau:V_r\to V_r$ be the shift operator: $\tau(f)(T)=f(qT)$.
Define 
$\alpha=T$ (operator of multiplication by $T$), $\gamma=-\tau^{-1}$ and $\beta=T^{-1}\circ(1-\tau)$.
One checks that all the relations of $A_q(S)$ are satisfied, and moreover, the seminorm on $A_q(S)$ induced
by the operator norm is multiplicative. (Similar considerations apply to the analytic quantum torus
derived from $\xi\eta=q\eta\xi$. Then the element $\sum_{n,m\in \Z}a_{nm}\xi^n\eta^m$
transforms the series $f=\sum_{i\in \Z}c_iT^i$ into $\sum_{n,m\in \Z}a_{nm}q^{nm}T^mf(q^nT)$). Rescaling the action
of $\alpha$ and $\gamma$ by arbitrary non-zero numbers one can adjust the action of $\beta$ in such a way
that the norms of operators $\alpha,\beta,\gamma$ ``cover" an open
neighborhood of the point $(1,1,1)$. More precisely, 
let us consider the map $f:M({\cal O}_q(S^{an}))\to \R^3$ defined by the formula
$f(x)=(a,b,c)$
where $a=\max(0,\log|\alpha|_x), b=\max (0, \log|\beta|_x),
c= \log|\gamma|_x=-\log|\alpha\beta-1|_x$. Here
$|\cdot|_x=\exp(-val_x(\cdot))$ denotes the mulitplicative seminorm
corresponding to the point $x\in M({\cal O}_q(S^{an}))$.
Then the image of $f$ is homeomorphic to $\R^2$, similarly to the case $q=1$ considered
in [KoSo1]. 
More precisely, let us decompose $M({\cal O}_q(S^{an}))=S_-\cup S_0\cup S_+$
 according to the sign of $\log |\gamma|_x$ where $x\in M({\cal O}_q(S^{an}))$. Then 
$$\begin{array}{lll}
f(S_-) & = & \{\,(a,b,c)\in {\R}^3\,|\,c<0,a\ge 0, b\ge 0, \,ab(a+b+c)=0\,\}\\
f(S_0) & = & \{\,(a,b,c)\in {\R}^3\,|\,c=0,a\ge 0, b\ge 0, \,ab=0\, \}  \\
f(S_+) & = &  \{\,(a,b,c)\in {\R}^3\,|\,c>0,a\ge 0, b\ge 0,\, ab=0\, \}
\end{array}$$
In fact the image of the map $f$ coincides with the image of the embedding
$j: {\R}^2\to {\R}^3$ given by  formula

$$j(x,y)=\left\{\begin{array}{lll}
(-x\,,\, \max(x+y,0)\,,\, -y\,) &
\mbox{ if } & x\le 0\\
 (\,0\,,\, x+\max(y,0)\,,\, -y\,) & \mbox{ if } & x\ge 0
\end{array}\right.$$

Proofs of the above observations are different from the case $q=1$. Indeed,
there are no one-dimensional modules over $A_q(S)$ corresponding
to the points of the surface $S$. Therefore it is not obvious that there are
multiplicative seminorms $x$ on $A_q(S)$ with the prescribed value of $f(x)$.
Seminorms on $A_q(S)$ arise from representations of this algebra in $k$-Banach
vector spaces: if $\rho:A_q(S)\to End_k(V)$ is such a representaion then we can
define $|a|_{\rho}=||\rho(a)||$, where $||\rho(a)||$ is the operator norm
in the Banach algebra $End_k(V)$  of bounded operators on $V$. Such seminorms 
are, in general, submultiplicative: $|ab|_{\rho}\le |a|_{\rho}|b|_{\rho}$.
We are interested in those which are multiplicative. This can be achieved, e.g.
by mapping of $A_q(S)$ into an algebra which  admits multiplicative seminorms.
We discussed above the homomorphisms $g_i$ of $A_q(S)$ into analytic quantum tori. 
Let us consider a different example of such  homomorphism.
Let $\delta=(\alpha\beta-1)\gamma$. One checks that $\delta$ is a central element
in the quantum affinoid algebra $B_q(S)$ generated by the first three relations 
for $A_q(S)$ (i.e. we drop the relation $\delta=1$).
Let us consider the quantum affinoid algebra $B$ generated by $\beta^{\pm 1},\gamma^{\pm 1}, \delta$
subject to the relations:

$$\beta\delta=\delta\beta, \gamma\delta=\delta\gamma, \gamma\beta=q\beta\gamma,$$
and such that $\beta^{-1}$ is inverse to $\beta$ and $\gamma^{-1}$ is inverse to
$\gamma$. There is an embedding of algebras $A_q(S)\to B/(\delta-1)B$ 
induced by the linear map $A_q(S)\to B$ such that $\beta$ and $\gamma$ are
mapped into the corresponding elements of $B$ and $\alpha\mapsto (1+\delta\gamma^{-1})\beta^{-1}$.
Notice that for any $r_1>0$, $r_2>0$, $r_3\ge 0$ one can define a multiplicative
norm $|\bullet|_{r_1,r_2,r_3}$ on $B$ such that
$|\sum_{n\in \Z,m\in \Z,l\in \Z_+}c_{nml}\beta^n\gamma^m\delta^l|_{r_1,r_2,r_3}=
max_{n,m,l}|c_{nml}|r_1^nr_2^mr_3^l$. Moreover, we can complete $B$ with respect
to this multiplicative norm and obtain a quantum affinoid algebra $B_{r_1,r_2,r_3}$. 
We can also invert $\delta$ and do the same construction. In this way we obtain
the quantum affinoid algebra denoted by $B_{r_1,r_2,r_3}^{(1)}$. Since $A_q(S)$ is embedded into 
the quotient of any of these
algebras by the central ideal, we obtain plenty multiplicative norms on $A_q(S)$ and 
on its completions.

We denote by $p$ the composition $j^{-1}\circ f$. In the case $q=1$ it is an analytic torus
fibration over the set $\R^2\setminus \{(0,0)\}$.
Let now $W$ be a small disc in $\R^2$ centered at the origin. 
We need to define the non-commutative
affinoid algebra ${\cal O}_q^{sing}(W)$. In the commutative case $q=1$ it is defined as
${\cal O}_{S^{an}}^{sing}(p_{\ast}^{-1}(W))=p_{\ast}({\cal O}_{S^{an}})(W)$. 

For each $i=1,2,3$ we define the $\C((t))$-affinoid algebras ${\cal O}_q^{sing}(p^{-1}(U_i))$
such that for every $x\in M({\cal O}_q^{sing}(p^{-1}(U_i))$ one has $p(x)\in U_i$ (it coincides with the intersection
of all completions of  $A_q(S)$ with respect to multiplicative seminorms $x$ such that 
$p(x)\in U_i,i=1,2,3$). Similarly we define algebras
${{\cal O}_q}^{sing}(p^{-1}(W))$ and ${{\cal O}_q}^{sing}(p^{-1}(W^0))$, where
$W^0=W\setminus \{(0,0)\}$ or, more generally, any ${\cal O}_q^{sing}(p^{-1}(U))$ for $U$ being an open subset
of $W$. Using homomorphisms $g_i, i=1,2,3$ one proves that if $U\subset U_i$ 
then ${\cal O}_q^{sing}(p^{-1}(U))$ is isomorphic to ${\cal O}_q^{mod}(\pi_i^{-1}(U))$. The latter is defined
as the set of series $\sum_{m,n\in \Z}c_{mn}\xi_i^m\eta_i^n$ such that $\pi_i(|\bullet|)\in U$
for any multiplicative seminorm $|\bullet|$ such that $sup_{m,n}(log|c_{mn}|+mlog|\xi_i|+nlog|\eta_i|<\infty)$, if $(x,y)\in U$.
The isomorphism of sheaves 
${{\cal O}_q^{sing}}|_{W\setminus \{(0,0)\}}\simeq {{\cal O}_q^{mod}}|_{W\setminus \{(0,0)\}}$ follows.
Details of this construction will be explained elsewhere.

\subsection{Trees, automorphisms and gluing}

As was explained in [KoSo1] in the case of $q=1$ and in [So1] in the case $|q|=1$,
one has to modify the canonical sheaf in order to glue it  with ${\cal O}_q^{sing}$.
Here we explain the construction following [KoSo1], [So1], leaving the details
to a separate publication.
The starting point for the construction is a subset ${\cal L}\subset B$ which is an infinite tree.
We called it {\it lines} in [KoSo1].

The definition is quite general. Here we discuss the $2$-dimensional case, while a much more complicated
higher-dimensional case was considered in the recent paper [GroSie1].
For a manifold $Y$
which carries a $\Z$-affine structure a {\it line $l$} is defined by
a continuous map $f_l: (0,+\infty)\to Y$ and a covariantly constant (with respect
to the connection which gives the affine structure)
nowhere vanishing integer-valued $1$-form
$\alpha_l\in \Gamma((0,+\infty),f_l^{\ast} ((T^\ast)^\Z)$.
A set ${\cal L}$ of lines is required to be decomposed into a disjoint union
${\cal L}={\cal L}_{in}\cup {\cal L}_{com}$ of {\it initial}
and {\it composite} lines. Each composite line is obtained as a result
of a finite number of ``collisions" of initial lines. A collision
is described by a $Y$-shape figure, where the bottom leg of $Y$ is a composite line,
while two other segments are ``parents" of the leg, so that the leg is obtained as a result
of the collision. A construction of the
set ${\cal L}$ satisfying the axioms from [KoSo1] was proposed in
[KoSo1], Section 9.3. Generalization to the higher-dimensional case can be found
in [GroSie1]. In the two-dimensional case the lines form an infinite tree embedded 
into $B$. The edges have rational slopes
with respect to the integral affine structure. The tree is dense in $B^0$.

With each line $l$ (i.e. edge of the tree) we  associate a continuous family of automorphisms
of stalks of sheaves of algebras
$\varphi_l(t): ({\cal O}^{can}_q)_{Y,f_l(t)}\to ({\cal O}^{can}_q)_{Y,f_l(t)}$.

Automorphisms $\varphi_l$ can be defined in the following way (see [KoSo1], Section 10.4).

First we choose affine coordinates in a neighborhhod of a point
$b\in B\setminus B^{sing}$, identifyin $b$ with the
point $(0,0)\in \R^2$.
Let $l=l_+\in {\cal L}_{in}$ be (in the standard affine coordinates)
a line in the half-plane $y>0$ emerging from $(0,0)$
(there is another such line $l_{-}$
in the half-plane $y<0$, see [KoSo1] for the details).
Assume that $t$
is sufficiently small. Then we define
$\varphi_l(t)$
on topological generators $\xi,\eta$ by the formula

$$\varphi_l(t)(\xi,\eta)=(\xi(1+\eta^{-1}),\eta).$$

In order to extend $\varphi_l(t)$ to the interval $(0,t_0)$, where
$t_0$ is not small, we cover the corresponding segment of $l$ by
open charts. Then a change of affine coordinates transforms
$\eta$ into a monomial multiplied by a constant from $(\C((t)))^{\times}$.
Moreover, one can choose the change of coordinates in such a way that
$\eta\mapsto C\eta$ where $C\in (\C((t)))^{\times}, |C|<1$ (such change
of coordinates preserve the $1$-form $dy$. Constant $C$ is equal to
$exp(-L)$, where $L$ is the length of the segment of $l$ between
two points in different coordinate charts).
Therefore $\eta$ extends analytically in a unique way to an element of
$\Gamma((0,+\infty), f_l^{\ast}(({\cal O}^{can}_q)^{\times}))$.
Moreover the norm $|\eta|$ strictly decreases as $t$ increases,
and remains strictly smaller than $1$. Similarly to [KoSo1], Section 10.4
one deduces that $\varphi_l(t)$ can be extended for all $t>0$.
This defines $\varphi_l(t)$ for $l\in {\cal L}_{in}$.

Next step is to extend $\varphi_l(t)$ to the case when $l\in {\cal L}_{com}$,
i.e. to the case when the line is obtained as a result of a collision
of two lines belonging to ${\cal L}_{in}$.
Following [KoSo1], Section 10, we introduce a group $G$ which contains
all the automorphisms $\varphi_l(t)$, and then prove the factorization theorem
(see [KoSo1], Theorem 6) which allows us to define $\varphi_l(0)$ in the case
when $l$ is obtained as a result of a collision of two lines $l_1$ and $l_2$.
Then we extend $\varphi_l(t)$ analytically for all $t>0$ similarly to the
case $l\in {\cal L}_{in}$.

More precisely, the construction of $G$ goes such as follows.
Let $(x_0,y_0)\in \R^2$ be a point,
$\alpha_1,\alpha_2\in (\Z^2)^\ast$ be $1$-covectors such that
$\alpha_1\wedge\alpha_2>0$.
Denote by $V=V_{(x_0,y_0),\alpha_1,\alpha_2}$ the closed angle
$$\{(x,y)\in \R^2|\langle \alpha_i, (x,y)-(x_0,y_0)\rangle\ge 0, i=1,2\,\}$$

Let ${\cal O}_q(V)$ be a $\C((t))$-algebra consisting
of series $f=\sum_{n,m \in \Z}c_{n,m}\xi^{n}\eta^{m}$,
such that $\xi\eta=q\eta\xi$ and $c_{n,m}\in \C((t))$ satisfy the condition that
for all $(x,y)\in V$ we have:

\begin{enumerate}
\item if $c_{n,m}\ne 0$ then $\langle (n,m),(x,y)-(x_0,y_0)\rangle\le 0$,
where we identified $(n,m)\in \Z^2$
with a covector in $(T_{p}^{\ast}Y)^{\Z}$;
\item $\log|c_{n,m}|+nx+my\to -\infty$ as long as $|n|+|m|\to +\infty$.
\end{enumerate}

For an integer covector $\mu=adx+bdy\in (\Z^2)^*$ we denote
by $R_{\mu}:=R_{(a,b)}$ the monomial $\xi^a\eta^b$. 
Then $R_{(a,b)}R_{(c,d)}=q^{ad-bc}R_{(c,d)}R_{(a,b)}=q^{-bc}R_{(a+c,b+d)}$.
We define a prounipotent group $G:=G(q,\alpha_1,\alpha_2,V)$ which consists of automorphisms of
${\cal O}_q(V)$  having the form $f\mapsto e^gfe^{-g}$ where

$$g=
\sum_{n_1,n_2\ge 0,n_1+n_2>0}c_{n_1,n_2}R_{\alpha_1}^{-n_1}R_{\alpha_2}^{-n_2} $$
where $c_{n_1,n_2}\in \C((t))$ and
$$\log|c_{n,m}|-n_1\langle \alpha_1, (x,y)\rangle-n_2\langle \alpha_2, (x,y)\rangle\le 0\,\,\,\,\forall \,(x,y)\in V$$
The latter condition is equivalent to
$\log|c_{n,m}|-\langle n_1 \alpha_1+n_2\alpha_2, (x_0,y_0)\rangle\le 0$.
The assumption $|q|=1$ ensures that the product is well-defined.

Let us consider automorphisms as above such that in the series for $g$ the ratio $\lambda=n_2/n_1\in [0,+\infty]_{\Q}:=\Q_{\ge 0}\cup\infty$ is fixed.  Such automorphism form a 
commutative subgroup $G_{\lambda}:=G_{\lambda}(q,\alpha_1,\alpha_2,V)\subset G$.
There is a natural map $\prod_{\lambda}G_{\lambda}\to G$, defined as in [KoSo1],
Section 10.2. The  factorization theorem proved in the loc. cit states that this
map is a bijection of sets.

\begin{exa} Let us consider the automorphism, discussed above:
$$\varphi(\xi,\eta)=(\xi(1+\eta^{-1}),\eta).$$
One can check that the transformation $\xi\mapsto \xi(1+\eta^{-1})$
has the form 
$$exp(Li_{2,q}(\eta^{-1})/(q-1))\xi exp(-Li_{2,q}(\eta^{-1})/(q-1)),$$
where  $Li_{2,q}(x)$ is the quantum dilogarithm function (see e.g. [BR]).
It satisfies the property
$(x;q)_{\infty}=exp(Li_{2,q}(-x)/(q-1))$, where $(a;q)_{N}=\prod_{0\le n\le N}(1-aq^n)$
for $1\le N\le \infty$. Using the formula
$(x;q)_{\infty}=\sum_{n\ge 0}{(-1)^nq^{n(n-1)/2}x^n\over{(q;q)_n}}$ one can show that
$lim_{q\to 1}Li_{2,q}(x)=Li_2(x)=\sum_{n\ge 1}(-1)^nx^n/n^2$, which is the ordinary dilogarithm
function (the latter appeared in [KoSo1], Section 10.4 in  the reconstruction problem
of rigid analytic K3 surfaces). The reader should notice that the quantum dilogarthm does not
behave well in the case $|q|=1$. Nevertheless one can define the corresponding group elements.
More details on that will be given elsewhere.
\end{exa}

Let us now assume that lines $l_1$ and $l_2$ collide at
$p=f_{l_1}(t_1)=f_{l_2}(t_2)$,
generating the line $l\in {\cal L}_{com}$.
Then $\varphi_l(0)$ is defined with the help of factorization theorem.
More precisely, we set $\alpha_i:=\alpha_{l_i}(t_i),\,\,i=1,2$ and the angle
$V$ is the intersection of certain half-planes $P_{l_1,t_1}\cap P_{l_2,t_2}$
defined in [KoSo1], Section 10.3. The half-plane $P_{l,t}$ is contained
in the region of convergence of $\varphi_l(t)$.
By construction, the elements
$g_0:=\varphi_{l_1}(t_1)$  and $g_{+\infty}:=\varphi_{l_2}(t_2)$ belong respectively
to $G_0$ and $G_{+\infty}$.
The we have:
$$g_{+\infty} g_0= {\textstyle \prod_\to}\left(
(g_\lambda)_{\lambda\in [0,+\infty]_\Q}\right)=
g_0\dots g_{1/2}\dots g_1 \dots
g_{+\infty}.$$

Each term $g_\lambda$ with $0<\lambda= n_1/n_2 <+\infty$
corresponds to the newborn line $l$ with the
direction covector
$n_1\alpha_{l_1}(t_1)+n_2\alpha_{l_2}(t_2)$.
Then we set  $\varphi_l(0):=g_{\lambda}$.
This transformation is defined by a series
which is convergent in a neighborhood of $p$, and using the analytic continuation
we obtain $\varphi_l(t)$ for $t>0$, as we said above. Recall that
every line carries an integer $1$-form $\alpha_l=adx+bdy$. By construction,
$\varphi_l(t)\in G_{\lambda}$, where $\lambda$ is the slope of $\alpha_l$.

Having automorphisms $\varphi_l$ assigned to lines $l\in {\cal L}$ we
proceed as in [KoSo1], Section 11, modifying the sheaf
${\cal O}^{can}_q$ along each line. We denote the resulting
sheaf by ${\cal O}^{mod}_q$. It is isomorphic to the
previously constructed sheaf
${\cal O}^{mod}_q$ in a neighborhood of the point $(0,0)$.

\begin{rmk}
The appearance of the dilogarithm function in the above example can be illustrated
in the picture of collision of two lines (say, $(x,0)$ and $(0,y)$ $ x,y\ge 0$) which
leads to the appearance of the new line, which is the diagonal $(x,x), x\ge 0$. Then the factorization
theorem gives rise to the five-term identity $g_{\infty}g_0=g_0g_1g_{\infty}$, which is the quantum version of the
famous five-term
identity for the dilogarithm function.

\end{rmk}

\section{$p$-adic quantum groups}

\subsection{How to quantize $p$-adic groups}

Let $L$ be a finite algebraic extension of the field ${\Q}_p$ of $p$-adic numbers,
and $K$ be a discretely valued subfield of the field
of complex $p$-adic numbers ${\C}_p$ containing $L$. All the fields carry non-archimedean norms,
which we will denote simply by $|\bullet|$ (sometimes we will be more specific, using the
notation like $|x|_K$ in order to specify which field we consider).
 We denote by ${\cal O}_L$ the ring
of integers of $L$ and by $m_L$ the maximal ideal of ${\cal O}_L$. Let $G$ be a locally
$L$-analytic
group, which is the group of $L$-points of a split reductive algebraic group ${\bf G}$ over $L$.
 Let $H\subset G$ be an open maximal compact subgroup. We would like to define quantum analogs of
the algebras $C^{la}(G,K), C^{la}(H,K)$ of locally analytic functions on $G$ and $H$, as well
as their strong duals $D^{la}(G,K)$, $D^{la}(H,K)$, which are the algebras of locally analytic distributions on $G$ and $H$ respectively
(see [SchT1], [Em1]). Modules over the algebras of locally analytic distributions were used in [SchT1-5], [Em1] for a description of
locally analytic admissible representations of locally $L$-analytic groups.
Our aim is to derive ``quantum" analogs of those results. In this paper we will
discuss definitions of the algebras only.

First of all we   are going to define locally analytic functions
and locally analytic distributions on the quantized
compact group $H$.
Let us explain our approach in the
case of the ``classical" (i.e. non-quantum) group $H$. We are going to present definitions
of $C^{la}(H,K)$ and $D^{la}(H,K)$ in such a way that
they can be generalized to the case of quantum groups. The difficulty which one
needs to overcome is to define everything using only two Hopf algebras: the universal enveloping
algebra $U(g), g=Lie(G)$ and the algebra $K[G]$ of {\it regular} functions on the algebraic group $G$.
Our construction consists of three steps.

1) For a ``sufficiently small" open compact subgroup $H_r\subset H$ we define the algebra
${C}^{an}(H_r,K)$
of {\it analytic} functions on $H_r$. Here $0<r\le 1$ is a parameter, such
that $H_1=H$, and if $r_1<r_2$ then $H_{r_1}\subset H_{r_2}$.
The strong dual to ${C}^{an}(H_r,K)$ is denoted by
$D^{an}(H_r,K)$. It is, by definition, the algebra of {\it analytic} distributions on $H_r$.
It can be described (see [Em1], Section 5.2) as a certain completion of the universal
enveloping algebra $U(g)$ of the Lie algebra $g=Lie(G)$.

2) For any $r\le 1$  we define a norm on the
algebra of regular functions $K[G]$ such that the completion with respect to this norm
is the algebra of {\it continuous} functions $C(H_r,K)$ on the group $H_r$.

3) In order to define {\it locally analytic} functions on $H$ we consider
a family of seminorms $|f|_l$ on $K[G]$, where $f\in K[G]$ and $l$ runs
through the set $D^{an}(H_r,K)$. More precisely, for every $l\in D^{an}(H_r,K)$ we define a seminorm $|f|_l=||(id\otimes l)(\delta(f))||$,
where $\delta$ is the coproduct on the Hopf algebra $K[G]$ and $||\bullet||$ is the norm
defined above on the step 2). The completion of $K[G]$ with respect to the topology defined
by the family of seminorms $|\bullet|_l$ is the algebra $C^{la}(H,K)$ of locally
analytic functions on $H$ defined in [SchT1]. The strong dual $D^{la}(H,K):=C^{la}(H,K)^{\prime}_b$ is the algebra of locally analytic distributions introduced in [SchT1].

We recall that the locally analytic representation theory of $G$ developed in [SchT1-6]
is based on the notion of coadmissible module over the algebra $D^{la}(H,K)$, where $H\subset G$
is an open compact subgroup. Therefore, from the point of view of representation theory, it suffices to quantize $D^{la}(H,K)$.

\subsection{Quantization of ``small" compact subgroups}

We would like to quantize $D^{la}(H,K)$ following the above considerations.
We will do that for the class of algebraic quantum groups
introduced by Lusztig (see [Lu1], [Lu2]).
Let us fix $q\in L^{\times}$ such that $|q|=1$ (this restriction is not necessary for algebraic quantization,
but it will be important when we discuss convergent series). We will assume
that there is $h\in {\cal O}_L$ such that $|h|<1, exp(h)=q$.
Let ${\bf G}$ be a semisimple simply-connected algebraic group over $\Z$,
associated with a  Cartan matrix $((a_{ij}))$ (more precisely, in order to be consistent
with the terminology of [Lu1] we  start with a {\it root datum} of finite type associated with a Cartan datum,
see [Lu1], Chapter 2. These data give rise to the Cartan matrix in the ordinary sense).
The algebraic group ${\bf G}(\C)$ of $\C$-points of ${\bf G}$ was quantized by Drinfeld (see e.g. [KorSo] Chapters 1,2).
We will need  a $\Z$-form of the
quantized algebraic group ${\bf G}$ introduced by Lusztig (see [Lu1]).
It allows us to define the quantized group over an arbitrary field. We need to be more specific when
speaking about ``quantized" group.
More precisely, following [Lu1] one can define  Hopf $L$-algebras $U_q(g_L)$ and $L[G]_q$,
which are the quantized enveloping algebra of the Lie algebra ${g_L}$ of ${\bf G}(L)$ and the algebra
of {\it regular} functions on the algebraic quantum group ${\bf G}(L)$ respectively.
Extending scalars to $K$ we obtain Hopf $K$-algebras $U_q(g_K)=K\otimes_L U_q(g_L)$ and
$K[G]_q=K\otimes_LL[G]_q$.
We will also need  $\Z$-forms of the above Hopf algebras, which will be denoted by $U:=U_A$ and $A[G]_v$ respectively.
The latter are Hopf algebras over the ring $A=\Z[v,v^{-1}]$, where $v$ is a variable. The algebras
$U_L$ and $L[G]_q$ are obtained by tensoring of $U$ and $A[G]_v$ respectively with $L$
in such a way that
$v$ acts on $L$ by multiplication by $q$.

As an $A$-module the algebra $U$ is isomorphic to the tensor product
$U\simeq U^+\otimes U^0\otimes U^-$ where $U^{\pm}$ are the quantized Borel subalgebras
and $U^0$ is the quantized Cartan subalgebra (see [Lu2]). Recall that $U^{+}$ (resp. $U^{-}$)
is an $A$-algebra generated by the divided powers $E_i^{(N)}$ (resp $F_i^{(N)}$) of the
Chevalley generators $E_i$ (resp. $F_i$) of the quantized enveloping
$\Q(v)$-algebra ${\bf U}$, where $1\le i\le n:=rank\, g_{L}$ (see [Lu2]). The algebra $U^0$ is generated over $A$ by the
generators
$K_i^{s_i}\binom{K_i}{n_i}_v$, where $\binom{K_i}{n_i}_v=
\prod_{1\le m\le n_i}{K_iv^{d_i(-m+1)}-K_i^{-1}v^{d_i(m-1)}\over{v^{d_im}-v^{-d_im}}}$.
Here $n_i\in \Z_+$, $s_i\in \{0,1\}$, and $K_i, 1\le i\le n$ are the standard Chevalley generators of ${\bf U}$. The integer numbers
$d_i\in \{1,2,3\}$ satisfy the condition that $((d_ia_{ij}))$ is a symmetric positive definite
matrix with $a_{ii}=2$, and $a_{ij}\le 0$ if $i\ne j$.
Recall that there is a canonical reduction of $U_A$ at $v=1$, which is a Hopf $\Z$-algebra
$U_{\Z}$. It is the universal enveloping of the integer Lie algebra $g_{\Z}$ of the corresponding
Chevalley group.
We will denote by $t_i\in g_{L}$ the generators at $v=1$ corresponding
to $K_i$, keeping the same notation $E_i^{(N)}, F_i^{(N)}$ for the rest of the generators of $g_{\Z}$.
Thus we have the standard decompostion $g_{\Z}=g_{\Z}^{+}\oplus h_{\Z}\oplus g_{\Z}^{-}$,
where the Lie algebra $g_{\Z}^{+}$ is generated (as a Lie algebra over $\Z$) by $E_i^{(N)}$, the Lie algebra $g_{\Z}^{-}$
is generated by $F_i^{(N)}$
and the commutative Lie algebra $h_{\Z}$ is generated by $\binom{t_i}{N}:=t_i(t_i-1)...(t_i-N+1)/N!$, $1\le i\le n$ (see [St], Theorem 2).
Lie algebra $g_L$ is generated by the standard Chevalley generators $E_i=E_i^{(1)},F_i=F_i^{(1)},t_i, 1\le i\le n$.
The Hopf algebra $U/(v-1)U$ is the universal enveloping algebra $U(g_L)$ of $g_L$.

In what follows, while keeping the above notation, we will assume for simplicity that $L=\Q_p$.

We will need the following extension of $U_q(g_L)$. Let us fix a basis $\{\alpha_i\}_{1\le i\le n}$
of simple roots of $g_L$, as well as invariant bilinear form on this Lie algebra such that $(\alpha_i,\alpha_j)=d_ia_{ij}$.

Let $h_{{\cal O}_L}=\oplus_{1\le i\le n}\Z_pt_i$. We  fix a global chart
$\psi: h_{{\cal O}_L}\to T^0$, where $T^0={\bf T}({\cal O}_L)$ is the maximal compact torus.
Then any element $a\in T^0$ can be written as an analytic function
$t=\psi(\sum_{1\le i\le n}x_it_i):=t(x_1,...,x_n)$, where $x_i\in \Z_p$.
Let us introduce a unital topological Hopf $L$-algebra ${U}_q^{an}(g_L)$ which is a Hopf  $L$-algebra generated by $E_i^{(N)}, F_i^{(N)},
 1\le i\le n, N\ge 1$ and the elements $t(x)=t(x_1,...,x_n)\in T^0$ as above, such that the relations between
$E_i^{(N)}, F_i^{(N)}$ are the same as
in $U_q(g_L)$, and
$t(x)E_i=E_it(x+v_i),
t(x)F_i=F_it(x-v_i),$
where $v_i=(a_{1i},a_{2i},...,a_{ni})$.
The elements $K_i^{\pm 1}=exp(\pm d_i h)$ (recall that $exp(h)=q$) belong to this
algebra and together with $E_i^{(N)}, F_i^{(N)}, 1\le i\le n, N\ge 1$ generate the Hopf algebra
isomorphic to ${U}_q(g_L)$.

There is a natural non-degenerate pairing
${U}_q^{an}(g_L)\otimes L[G]_q\to L$ which extends the natural non-degenerate pairing $U_A\otimes A[G]_v\to A$
defined in [Lu2].
Extending scalars we obtain the algebra ${U}_q^{an}(g_K)$  and the pairing
${U}_q^{an}(g_K)\otimes K[G]_q\to K$.

{\it For the rest of this subsection we will assume that $d_i=1$, i.e. $((a_{ij}))$
is symmetric, and $L=\Q_p$. These conditions can be relaxed. We make them in order
to simplify formulas}.

There is  a natural action  $U_q^{an}(g_K)\otimes K[G]_q\to K[G]_q$ (right  action) given by the formula
$l(f)=(id\otimes l)(\delta(f))$, where $l\in U_q^{an}(g_K), f\in K[G]_q$ and
$\delta: K[G]_q\to K[G]_q\otimes K[G]_q$ is the coproduct.

Recall that $K[G]_q\simeq 
\oplus_{\Lambda}m_{\Lambda}V(\Lambda)$ which is the sum of irreducible finite-dimensional highest weight 
$U_q(g_K)$-modules $V(\Lambda)$ with multiplicities $m_{\Lambda}$. This is also an isomorphism of $U_q^{an}(g_K)$-modules.
Each element $E_i^{(N)}, F_i^{(N)}$ acts locally nilpotently on $K[G]_q$, 
while each $t(x)$ acts as a  semi-simple linear map.

Let $R=\oplus_{1\le i\le n}\Z\alpha_i$ be the set of roots of $g_{\Z}$. 
We denote by $R^+$ (resp $R^-$) the set of positive (resp. negative) roots.
We will often write $\alpha>0$ (resp. $\alpha<0$) if $\alpha\in R^+$ (resp. $\alpha\in R^-$).
Following [KorSo], Chapter 4, or [Lu1], Chapter 3, 41, 
one can construct quantum root vectors $E_{\alpha}^{(N)}, F_{\alpha}^{(N)}\in U_q(g_K), \alpha>0, N\ge 1$,
such that $E_{\alpha_i}^{(N)}=E_i^{(N)}, F_{\alpha_i}^{(N)}=F_i^{(N)}$ (in order to keep
track of integrality of the coefficients we are going to use the formulas from [Lu1]).
Let us fix a convex linear order on the set of roots,
such that all negative roots preceed all positive roots (convexity
means that  $\alpha<\alpha+\beta<\beta$ for positive roots and the oppoiste inequalities for negative roots) .

For every $0<r\le 1$ we  define $U_q(g_K)(r)$ as a $K$-vector space consisting
of series 
$$\xi=\sum_{m\in \Z^n, \alpha>0, s_{\alpha},p_{\alpha}\ge 0}c_{m,s,p}t^m/m!\prod_{\alpha>0}F_{\alpha}^{(p_{\alpha})}E_{\alpha}^{(s_{\alpha})},$$
such that $c_{m,s,p}\in K, t^m/m!=t_1^{m_1}/m_1!...t_n^{m_n}/m_n!, |c_{m,s,p}|r^{-(|s|+|p|+|m|)}\to 0$ as
$|m|+|s|+|p|\to \infty$. Here and below $m,s,p$ denote multi-indices.
We define $|\xi|_r=sup_{m,s,p}|c_{m,s,p}|r^{-(|s|+|p|+|m|)}$.
Let $t_{\alpha_i}(x), 1\le i\le n$ be an ordered basis of $T^0$ (see [SchT1], Section 4).
Then, as a topological $K$-vector space (with the topology defined by the norm $|\bullet|_r$) the space $U_q(g_K)(r)$
is isomorphic to the $K$-vector space of infinite series 
$$\eta=\sum_{N_i,M_i\ge 0, l_i\in \Z}b_{M,N}\prod_{1\le i\le n}(t_{\alpha_i}(x)-1)^{l_i}F_i^{(M_i)}E_i^{(N_i)},$$
such that $M=(M_i), N=(N_i)$, and $|b_{M,N,l}||t_{\alpha_i}(x)-1|_rr^{-(|M|+|N|)}\to 0$ as $|M|+|N|+|l|\to \infty$, equipped with the
norm defined by 
$$|\eta|^{\prime}_r=sup_{M,N,l}|b_{M,N,l}||t_{\alpha_i}(x)-1|_rr^{-(|M|+|N|)}. $$
It is easy to see that $U_q(g_K)(r)$ is a $K$-Banach vector space.
It contains Banach vector subspaces $U_q^+(g_K)(r)$ (resp. $U_q^-(g_K)(r)$)
which are closures of vector subspaces generated by all the elements $E_{\alpha}^{(N)}$ (resp. $F_{\alpha}^{(N)}$). 
It also contains an analytic neighborhood of $1\in T^0$, which is an analytic group
isomorphic to the ball of radius $r$ in the Lie algebra $h_{{\cal O}_L}$. The latter is an analytic 
Lie group via Campbell-Hausdorff formula. We can always assume that $r$ belongs to the algebraic
closure  of $L$, thus the corresponding analytic groups are in fact affinoid.

\begin{prp}
The norm $|\xi|_r$ (equivalently the norm $|\eta|^{\prime}_r$) gives rise  
to a Banach $K$-algebra structure on $U_q(g_K)(r)$. 
\end{prp}

Similarly to the case $q=1$ (see [Em1], Section 5.2) one can ask whether the algebra
$U_q(g_K)(r)$ corresponds to a ``good" analytic group. Let us consider the completion
of the tensor product $U_q(g_K)(r)\otimes U_q(g_K)(r)$ with respect to the minimal Banach norm.
Then we have the following result.

\begin{prp} The Hopf algebra structure on $U_q^{an}(g_K)$ admits a continuous extension to $U_q(g_K)(r)$,
making it into a topologocal Hopf algebra.

\end{prp}

Let us consider the topological $K$-algebra $U_q^{(1)}(g_K)$ which is the projective limit of
$U_q(g_K)(r)$ for all $0<r< 1$. Then we have the following result, which is a corollary of the
previous Proposition.

\begin{prp} The Hopf algebra structure on $U_q^{an}(g_K)$ admits a continuous extension to $U_q^{(1)}(g_K)$,
making it into a topological Hopf algebra.

\end{prp}

Since the elements $E_{\alpha}, F_{\alpha}$
act locally nilpotently on $K[G]_q$, there is a well-defined action
of $U_q(g_K)(r)$ on $K[G]_q$, which extends to the action of $U_q^{(1)}(g_K)$ on $K[G]_q$.
Notice that the pairing $U_q(g_K)(r)\otimes K[G]_q\to K, (l,f)\mapsto l(f)$ is non-degenerate.
In particular we can define the norm on $K[G]_q$ by the formula
$||f||_r=sup_{l\ne 0}{|l(f)|\over {|l|_r}}, l\in U_q(g_K)(r)$.

Let now $H_r, r=p^{-N}$ be a ``small" compact open subgroup of $G$. This means that
the exponential map $exp: \Z^d=h_{\Z}\oplus g_{\Z}^+\oplus g_{\Z}^-\to G$
defines an analytic isomorphism
$B(0,r)\to H_r$, where $B(0,r)\subset \Z_p^d$ is the ball consisting of points 
$(x_i,x_{\alpha},y_{\alpha})\in \Z_p^d, 1\le i\le n, \alpha>0$ such that
$x_{\alpha},y_{\alpha}\in p^N\Z_p,x_i\in p^N\Z_p$, for all $\alpha>0, 1\le i\le n$.

\begin{defn} The space of analytic functions on the quantum group $H_r$
(notation ${C}^{an}(H_r,K)_q$) is the completion of $K[G]_q$ with respect to the
norm $||\bullet||_r$.
\end{defn}

\begin{prp} The space ${C}^{an}(H_r,K)_q$ is  a Banach Hopf $K$-algebra.

\end{prp}

\begin{defn} The algebra of analytic distributions on the quantum group $H_r$
(notation $D^{an}(H_r,K)_q$) is the strong dual to ${C}^{an}(H_r,K)_q$.

\end{defn}

One can define a norm $||\bullet||$ on $K[G]_q$ such that the completion
with respect to this norm is by definition the algebra $C(H,K)_q$ of continuous functions
on the open maximal compact subgroup $H=H_1$ (in the case of $q=1$ this is a
theorem, not a definition). Then we proceed as follows.

Any linear functional $l\in D^{an}(H_r,K)_q$ defines a seminorm $|\bullet|_l$ on
$K[G]_q$ such that
$$|f|_l=||(id\otimes l)\delta(f)||.$$
The collection of seminorms $|\bullet|_l, l\in D^{an}(H_r,K)_q$ gives rise to
a locally convex topology on $K[G]_q$.

\begin{defn} The space ${C}^{an}(H, H_r, K)_q$ of functions on
the quantum group $H$ which are locally
analytic with respect to the quantum group $H_r$ is  the completion
of $K[G]_q$ in the topology defined by the collection of seminorms
$|\bullet|_l$.

\end{defn}

\begin{defn} a) The space $C^{la}(H,K)_q$ of locally analytic functions on the quantum group $H$ is the inductive limit $\varinjlim_{r\le 1}{\cal O}^{an}(H, H_r, K)$
(i.e. it consists of functions on quantum group $H$ which are locally analytic with respect to some $H_r, r<1)$.

b) The space $D^{la}(H,K)_q$ of locally analytic distributions on the quantum group $H$
is  the strong dual to $C^{la}(H,K)_q$.

\end{defn}

Since some details related to the proof of the following results are not finished, I formulate
it as a conjecture.

\begin{conj} Both spaces $C^{la}(H,K)_q$ and $D^{la}(H,K)_q$ are topological Hopf
$K$-algebras. Furthermore, $D^{la}(H,K)_q$ is a Frech\'et-Stein algebra in the sense
of [SchT1].

\end{conj}

In the next subsection we will explain the definition of the norm $||\bullet||$ in the case
of the group $SL_2(\Z_p)$. The general case is similar, but requires more details.
It will be considered in [So2].

\subsection{The $GL_2(\Z_p)$-case }

We will use the notation
$K[GL_2(\Q_p)]_q$ for the algebra of regular $K$-valued functions
on the algebraic quantum group $GL_2(\Q_p)$. It is known 
(see  [KorSo], Chapter 3) that
$K[GL_2(\Q_p)]_q$ is generated by generators $t_{ij}, 1\le i,j\le 2$ subject
to the relations
\begin{eqnarray}
\label{1.1-1.2}
t_{11}t_{12}=q^{-1}t_{12}t_{11}, & t_{11}t_{21}=q^{-1}t_{21}t_{11}, \\
\label{1.3-1.4}
t_{12}t_{22}=q^{-1}t_{22}t_{12}, & t_{21}t_{22}=q^{-1}t_{22}t_{21}, \\
\label{1.5-1.6}
t_{12}t_{21}=t_{21}t_{12}, & t_{11}t_{22}-t_{22}t_{11}=
\left(q^{-1}-q\right)t_{12}t_{21},
\end{eqnarray}

The element $det_q=t_{11}t_{22}-q^{-1}t_{12}t_{21}$ generates
the center of the above algebra. As a result, the algebra
$K[SL_2(\Q_p)]_q$ of regular functions on quantum group $SL_2(\Q_p)$
is obtained from the above algebra by adding one more equation

\begin{equation}
\label{III-1.7eq}
t_{11}t_{22}-q^{-1}t_{12}t_{21}=1.
\end{equation}

We are going to use the ideas of the representation theory
of quantized algebras of functions (see [KorSo]).

Let $V$ be a separable $K$-Banach vector space. This means that
$V$ contains a dense $K$-vector subspace spanned by the orthonormal
basis $e_m, m\ge 0$ (orthonormal means that $||e_m||=1$ for all $m$).
Let us consider the following representations $V_c, c\in K$ of $K[GL_2(\Q_p)]_q$
in $V$ (cf. [KorSo], Chapter 4, Section 4.1):

$$t_{11}(e_m)=a_{11}(m)e_{m-1}, t_{21}(e_m)=a_{21}(m)e_m,$$

$$t_{12}(e_m)=a_{12}(m)e_m, t_{22}(e_m)=a_{22}(m)e_{m+1},$$

$$det_q=c.$$

Here $a_{ij}(m)\in K$ and $e_m=0$ for $m<0$. In particular
the line $Ke_0$ is invariant with respect to the subalgebra
$A_+$ generated by $t_{11}, t_{21}$. Let us assume that not all
$a_{11}(m)$ are equal to zero. Then
it is easy to see from the commutation relations between $t_{ij}$ that
$$ a_{21}(m)=a_{21}(0)q^{-m},a_{12}(m)=a_{12}(0)q^{-m}, m\ge 0.$$ 

Moreover 
$$a_{11}(m+1)a_{22}(m)-a_{11}(m)a_{22}(m-1)=(q^{-1}-q)q^{-2m}h_0,$$
where $h_0=a_{21}(0)a_{12}(0)$. 
Let $s(m)=a_{11}(m)a_{22}(m-1), m\ge 1$. Then we have
$$s(m+1)-s(m)=(q^{-1}-q)q^{-2m}h_0, s(1)=(q^{-1}-q)h_0.$$
It follows that
$$s(m)=(q^{-1}-q)(1+q^{-2}+...+q^{-2(m-1)})h_0=q(q^{-2m}-1)h_0,$$
for all $m\ge 1$.
Since
the quantum determinant is equal to $c$,  we have
$$s(m+1)=a_{11}(m+1)a_{22}(m)=c+q^{-2m-1}h_0.$$
Comparing two formulas for $s(m)$ we see that
$$h_0=a_{21}(0)a_{12}(0)=-cq^{-1}.$$

{\bf From now on we will assume that $|1-q|<1$}.

Then the operators $t_{12}$ and $t_{21}$ are bounded.
We also have $|s(m)|=|c(q^{-m}-1)|=|c|, m\ge 1$. 

Assume that $a_{21}(0)\ne 0$.
Then the above representations (which are algebraically irreducible
as long as $q$ is not a root of $1$) depend on the parameters
$a_{21}(0), a_{11}(m), a_{22}(m), m\ge 0$ subject to the relations
$a_{11}(m)a_{22}(m-1)=c(1-q^{-2m})$. We will further specify 
restrictions on these parameters. The idea is the same as in
[KorSo], Chapter 3, where in order to define continuous
functions on the quantum group $SU(2)$ we  singled out 
irreducible representations
of $\C[SL_2(\C)]_q$ corresponding to the intersection of the
group $SU(2)$ with the big Bruhat cell for $SL_2(\C)$. 
This intersection is the union of symplectic leaves of the Poisson-Lie group 
$SU(2)$. Kernel of an irreducible representation defines a symplectic leaf (``orbit method") which explains the relationship 
of representation theory and symplectic geometry.
Notice that in the case $q=1$ one can define the algebra
of continuous function $C(SU(2))$ in the following way.
For any function $f\in\C[SL_2(\C)]$ one takes its restriction
to the above-mentioned union of symplectic leaves.
Since the latter is dense in $SU(2)$, the completion of the algebra
$C[SL_2(\C)]$ with respect to the sup-norm taken over all irreducible representations corresponding to the symplectic leaves
is exactly $C(SU(2))$.
Now we observe that symplectic leaves in $SL_2(\C)$
are  algebraic subvarieties, therefore they exist over any field.
We will use the same formulas in the case of  any $p$-adic field $L$ (in this section
we take $L=\Q_p$)..
In order to specify a symplectic leaf in $GL_2(\C)$
we need in addition to fix the value of the determinant (it belongs
to the center of the Poisson algebra $\C[GL_2(\C)]$).

Let us recall (see [KorSo], Chapter 3) that to every element
$t,c\in K^{\times}$ one can assign a $1$-dimensional representation  $W_{c,t}=\Q_pe_0$
of $K[GL_2(L)]_q$ such that $t_{11}(e_0)=te_0, t_{22}(e_0)=ct^{-1}e_0$,
and the rest of generators act on $e_0$ trivially.
Recall (see [KorSo], Chapter 1) that complex $2$-dimensional symplectic leaves of
$GL_2(\C)$ are algebraic subvarieties $S_{c,t}$ given by the equations:
$$t_{11}t_{22}-t_{21}t_{12}=c, t_{12}=t^2t_{21},$$
where $c,t$ are non-zero complex numbers.
We define symplectic leaves over $\Q_p$ by the same formulas,
taking $c,t\in \Q_p$.
In order to define the norm of the restriction of a regular
function $f\in \C[GL_2(\Q_p)]$ on $GL_2(\Z_p)$ we can choose
a subset in the set of symplectic leaves $S_{c,t}$ such that
the union of their intersection with $GL_2(\Z_p)$ is dense
in the latter group. It suffices to take those leaves
$S_{c,t}$ for which $|c|\le 1, t\in \Z_p^{\times}$, and both $t_{12}$ and $t_{21}$ are non-zero. 

Let us consider infinite-dimensional represenation $V_{c,t}$ as above
for which $a_{12}(0)=t^2a_{21}(0), |c|\le 1$ for a fixed $i\ge 0$,
and $a_{21}(0)\ne 0$. We will also assume that
the norm of the operators corresponding to $t_{11}$ and $t_{22}$
is less or equal than $1$. It follows from the equality
$t^2a_{21}^2(0)=-cq^{-1}$ that $|a_{21}(0)|=|c|$, hence the norm
of the operators corresponding to $t_{12}$ and $t_{21}$
is less or equal than $|c|\le 1$. It follows that
the norm of the operator $\pi_{c,t}(f)$  corresponding to
an element $f\in K[GL_2(\Q_p)]_q$ acting in $V_{c,t}$
is bounded from above as $V_{c,t}$ run through the set of irreducible
representations with the above restrictions on $c,t$.
In addition, we are going to consider 
only those $c\in K^{\times}$ for which $-cq^{-1}$ is a square
in $K$.
We define the norm $||f||_{GL_2(\Z_p),q}$,
$f\in K[GL_2(\Q_p)]_q$ as the supremum of norms of the operators
$\pi_{c,t}(f)$  corresponding
an element $f$ in all representations $V_{c,t}$ as above. This is the desired sup-norm
which we  used in our definition of the algebra of locally analytic functions.

\vspace{3mm}

{\bf References}

\vspace{2mm}

[BR] V. Bazhanov, N. Reshetikhin, Remarks on the quantum dilogarithm, J. Phys. A: Math. Gen. 28, 1995. 2217-2226.

\vspace{2mm}

[Be1] V. Berkovich, Spectral theory and analytic geometry over non-archimedean fields, Amer. Math. Soc.,
1990.

\vspace{2mm}

[Be2] V. Berkovich, Smooth $p$-adic analytic spaces are locally contractible, Invent. Math., 137, 1999, 1-84.

\vspace{2mm}

[Be3] V. Berkovich, Smooth $p$-adic analytic spaces are locally contractible II, In: Geometric aspects of Dwork theory, Berlin, 293-370.

[Fro] H. Frommer, The locally analytic principal series of split reductive groups, preprint 265, Univerity of Muenster, 2003.

\vspace{2mm}

[GroSie1] M. Gross, B. Siebert, From real affine geometry to complex geometry, math.AG/0703822.
\vspace{2mm}

[Em1] M.Emerton,  Locally analytic vectors in representations of locally p-adic analytic groups,
http://www.math.northwestern.edu/~emerton/preprints.html, to appear in  Memoirs of the AMS.

\vspace{2mm}

[Em2] M.Emerton, Locally analytic representation theory of p-adic reductive groups: 
 A summary of some recent developments,
to appear in the proceedings of the LMS symposium on L-functions and Galois representations .

\vspace{2mm}
[KoSo1] M. Kontsevich, Y. Soibelman, Affine structures and non-archimedean analytic spaces,
math.AG/0406564.

\vspace{2mm}

[KoSo2] M. Kontsevich, Y. Soibelman, Homological mirror symmetry and torus fibrations,
math.SG/0011041.

\vspace{2mm}

[KorSo] L.Korogodsky, Y. Soibelman, Algebras of functions on quantum groups I.
AMS, 1998.

\vspace{2mm}

[Lu1] G. Lusztig, Introduction to quantum groups, Birkhauser, 1993.

\vspace{2mm}

[Lu2] G. Lusztig, Quantum groups at roots of 1, Geometriae Dedicata, v. 35, 89-114, 1989.

\vspace{2mm}

[Lu3] G. Lusztig, Finite-dimensional Hopf algebras arising from quantized universal enveloping algebras,
J. Amer. Math. Soc., 3:1, 1990, 257-296.

\vspace{2mm}

[Sch-NFA] P. Schneider, Non-archimedean functional analysis, Springer.
\vspace{2mm}

[SchT1] P. Schneider, J.Teitelbaum, Algebras of $p$-adic distributions and
admissible representations, math.NT/0206056.
 
 \vspace{2mm}

[SchT2] P. Schneider, J.Teitelbaum,p-adic Fourier Theory,math.NT/0102012.
 
\vspace{2mm}

[SchT3] P. Schneider, J.Teitelbaum, U(g)-finite locally analytic representations,
ath.NT/0005072.

\vspace{2mm}

[SchT4] P. Schneider, J.Teitelbaum, Banach space representations and Iwasawa theory,
math.NT/0005066.

\vspace{2mm}
     
[SchT5] P. Schneider, J.Teitelbaum, Locally analytic distributions and p-adic representation theory, with applications to $GL_2$, math.NT/9912073.
 
 \vspace{2mm}
 
[SchT6] P. Schneider, J.Teitelbaum, p-adic boundary values,  
math.NT/9901159.

\vspace{2mm}

[So1] Y. Soibelman, On non-commutative spaces over non-archimedean fields, math.QA/0606001.

\vspace{2mm}

[So2] Y. Soibelman, Quantum p-adic groups and their representations, in preparation.

\vspace{2mm}

[So3] Y. Soibelman, Algebras of functions on a compact quantum group and its representations,
Leningrad Math. Journ., vol.2, 1991, 193-225.
\vspace{2mm}

[SoVo] Y. Soibelman, V. Vologodsky, Non-commutative compactifications and elliptic curves,
math.AG/0205117, published in Int. Math.Res. Notes, 28 (2003).

\vspace{5mm}

\end{document}